\documentclass[preprint,12pt]{elsarticle}

\usepackage{amssymb}

\usepackage{setspace}
\newtheorem{thm}{Theorem}[section]
\newtheorem{cor}[thm]{Corollary}
\newtheorem{lem}[thm]{Lemma}
\newtheorem{prop}[thm]{Proposition}
\newtheorem{rem}[thm]{Remark}

\begin{document}

\begin{frontmatter}



\title{Positive solutions for a class of singular quasilinear Schr\"{o}dinger equations with critical Sobolev exponent}


\author{Zhouxin Li}
\ead{lzx@math.pku.edu.cn}
\address{Department of Mathematics and Statistics, Central South University,\\ Changsha 410083, PR China}


\begin{abstract}
In this paper we prove the existence of positive solutions of the following singular quasilinear Schr\"{o}dinger equations at critical growth
\begin{eqnarray*}
-\Delta u-\lambda c(x)u-\kappa\alpha(\Delta(|u|^{2\alpha}))|u|^{2\alpha-2}u
    =|u|^{q-2}u+|u|^{2^*-2}u,\quad u\in{D^{1,2}(\mathbb{R}}^N),
\end{eqnarray*}
via variational methods, where $\lambda\geq0$, $0<\alpha<1/2$, $2<q<2^*$. It is interesting that we do not need to add a weight function to control $|u|^{q-2}u$.
\end{abstract}


\end{frontmatter}


\section{ Introduction and main results}
\setcounter{equation}{0}
\label{sec:5-1}
In this paper, we consider the following quasilinear Schr\"{o}dinger equation introduced in \cite{LiWa03,LiWW03}
\begin{equation}\label{eq.sch5-1}
  i\partial_tz=-\Delta z+w(x)z-l(|z|^2)z-\kappa\Delta h(|z|^2)h'(|z|^2)z,\quad x\in{\mathbb{R}}^N,
\end{equation}
where $w(x)$ is a given potential, $\kappa>0$ is a constant, $N\geq3$.
$h,~l$ are real functions of essentially pure power form.

Eq.(\ref{eq.sch5-1}) comes from mathematical physics and was used to model some
physical phenomena.
If $\kappa=0$, Eq.(\ref{eq.sch5-1}) is a semilinear problem which has been extensively studied.
If $\kappa>0$, it is a quasilinear problem which has many applications in physics.
It is known that the case of $h(s)=s$ was used for the superfluid film equation in plasma physics by Kurihura in \cite{Kuri81}.
It also appears in plasma physics and fluid mechanics \cite{LiSe78}, in the theory of Heisenberg ferromagnetism and magnons \cite{KoIK90,QuCa82} in dissipative quantum mechanics \cite{Hass80} and in condensed matter theory \cite{MaFe84}.

If we consider solutions of the form
$z(x,t)=\exp(-iet)u(x)$, which are called standing waves,
we observe that this $z(x,t)$ satisfies Eq.(\ref{eq.sch5-1}) if and only
if the function $u(x)$ solves the equation
\begin{eqnarray}\label{eq.sch5-2}
-\Delta u+V(x)u-\kappa\alpha(\Delta(h(|u|^2)))h'(|u|)^{2}u
    =l(|u|^2)u,\quad x\in{\mathbb{R}}^N,
\end{eqnarray}
where $V(x)=w(x)-e$ is the new potential function.

In recent years, the case $h(s)=s$ has been extensively studied under different conditions on the potential $V(x)\geq0$ and the nonlinear perturbation $l(u)$, one can refer to \cite{LiWW03,CoJe04,OMS07,OMS10} and some references therein.
The difficulty of Eq.(\ref{eq.sch5-2}) lies in the unbounded operator. In order to overcome this difficulty, Liu and Wang etc. in \cite{LiWW03} defined a change of variable and change the problem to a semilinear one. More precisely, they used the change of variable $v=f^{-1}(u)$ with $f$ defined by ODE: $f'(t)=(1+2f^2(t))^{-1/2}, ~t\in(0,+\infty)$ and $f(t)=-f(-t), ~t\in(-\infty,0)$. Then they proved the existence of positive solutions in an corresponding Orlicz space.
This method was widely used in the studies of such kind of problems thereafter, for examples \cite{CoJe04,OMS07,OMS10}. But in the latter literatures the working space is the usual Sobolev space $H^1(\mathbb{R}^N)$.

It is also interesting to study Eq.(\ref{eq.sch5-2}) with the nonlinearity $l(s)$ is at critical growth. In \cite{LiWa03}, Liu and Wang pointed out that the number $2(2^*)$ behaves like critical exponent for Eq.(\ref{eq.sch5-2}).
In \cite{SiVi10}, Silva and Vieira proved the existence of solutions of Eq.(\ref{eq.sch5-2}) with a general $l(|u|^2)u=K(x)u^{2(2^*)-1}+g(x,u)$.

To our knowledge, there are few literatures study more general problem $h(s)=s^{\alpha}$ with $\alpha>1/2$. We mention that in \cite{LiWa03}, the existence results are obtained through a constrained minimization argument for $l(u)$ at subcritical growth. In \cite{Moam06}, Moameni consider the problem for $l(u)$ at critical growth under radially symmetric conditions. But note that such assumptions enable the study to avoid the difficulty of losing compactness caused by Sobolev imbedding.
In \cite{LiZh13}, Li and Zhang studied the problem that $h(s)=s^{\alpha}$,
$l(s)=s^{(q-2)/2}+s^{(2^*-2)/2}$, where $\alpha>1/2,~2(2\alpha)\leq q<2^*(2\alpha)$, $2^*=2N/(N-2)$, and proved the existence of positive solution. Here and in the following, we always denote $2(2\alpha)=2\times2\alpha$, $2^*(2\alpha)=2^*\times2\alpha$.

Compare to \cite{LiZh13}, in this paper, we are interested in the problem with $h(s)=s^{\alpha}$,
$l(s)=s^{(q-2)/2}+s^{(2^*-2)/2}$, where $0<\alpha<1/2,~2\leq q<2^*$.
It was used to models the self-channeling of high-power ultrashort laser in matter \cite{BoGa93}. 
Moreover, we consider problems that $V(x)<0$.

Denote
\[
X:= D^{1,2}({\mathbb{R}}^N)=\{u\in L^{2^*}({\mathbb{R}}^N):|\nabla u|\in L^2({\mathbb{R}}^N)\}
\]
equipped with the norm $\|u\|^2_X=\int_{{\mathbb{R}}^N}|\nabla u|^2\mbox{d}x$.

Let $\lambda c(x):=-V(x)>0,~\lambda\geq0$.
In this paper, we consider the following equation
\begin{eqnarray}\label{eq.sch5-3}
-\Delta u-\lambda c(x)u-\kappa\alpha(\Delta(|u|^{2\alpha}))|u|^{2\alpha-2}u
    =|u|^{q-2}u+|u|^{2^*-2}u,\quad x\in{\mathbb{R}}^N.
\end{eqnarray}

Note that when $0<\alpha<1/2$, the operator of second order is singular in the equation, this cause one of the main difficulty of the study.
Another difficulty of the study is caused by the nonlinear term $|u|^{q-2}u$ in the equation since $D^{1,2}(\mathbb{R}^N)$ does not imbed into $L^{q}(\mathbb{R}^N)$.

Let $f(t)=|u|^{q-2}u+|u|^{2^*-2}u$, $2\leq q<2^*$.
We want to find weak solutions to Eq.(\ref{eq.sch5-3}). By {\it weak solution}, we mean a function $u$ in $X$ satisfying that,
for all $\varphi\in C_0^{\infty}({\mathbb{R}}^N)$, there holds
\begin{eqnarray}\label{def.I'}
\int_{{\mathbb{R}}^N}\nabla u\nabla\varphi-\lambda\int_{{\mathbb{R}}^N}c(x)u\varphi+\kappa\alpha\int_{{\mathbb{R}}^N}\nabla(|u|^{2\alpha})\nabla(|u|^{2\alpha-2}u\varphi)
 =\int_{{\mathbb{R}}^N}f(u)\varphi.
\end{eqnarray}

According to the variational methods, the weak solutions of (\ref{eq.sch5-3}) corresponds the critical points
of the functional $I: X\to{\mathbb{R}}$ defined by
\begin{eqnarray}\label{def.I}
I(u)=\frac{1}{2}\int_{{\mathbb{R}}^N}(1+2\kappa\alpha^2|u|^{2(2\alpha-1)})|\nabla u|^2-\frac{\lambda}{2}\int_{{\mathbb{R}}^N}c(x)u^2
-\int_{{\mathbb{R}}^N}F(u),
\end{eqnarray}
where $F(t)=\int_0^uf(s)\mbox{d}s$. For $u\in X$, $I(u)$ is lower semicontinuous when $0<\alpha<1/2$, and not differentiable in all directions $\varphi\in X$. In order to use the classical critical point theorem, we will use a change of variable to reformulate functional $I$.

Let $\beta(t)=(1+2\kappa\alpha^2|t|^{2(2\alpha-1)})^{1/2}$, then $\beta(t)$ is monotone and decreasing in $t\in(0,+\infty)$. For $t_0>0$ sufficiently small,
we have
\[
\int_0^{t_0}\beta(s)ds\leq2\alpha\sqrt{\kappa}\int_0^{t_0}s^{2\alpha-1}ds=\sqrt{\kappa}t_0^{2\alpha}.
\]
This make it possible for us to define an invertible, odd, $C^1$ function $h:{\mathbb{R}}\to{\mathbb{R}}$ by
\begin{eqnarray}\label{def.h}
v=h^{-1}(u)=\int_0^u\beta(s)ds,
\end{eqnarray}
where $h^{-1}$ is the inverse function of $h$.
For simplicity of notation we may assume that $\kappa\alpha=1$.

Inserting $u=h(v)$ into (\ref{def.I}), we get another functional $J$ defined on $X$ given by
\begin{eqnarray}\label{def.J}
J(v)=I(h(v))=\frac{1}{2}\int_{{\mathbb{R}}^N}|\nabla v|^2
-\frac{\lambda}{2}\int_{{\mathbb{R}}^N}c(x)h(v)^2-\int_{{\mathbb{R}}^N}F(h(v)).
\end{eqnarray}
In Sect.\ref{sec:5-2} (see Proposition \ref{prop.J}) we prove that $J$ is well defined on $X$, and is continuous in $X$. Moreover, it is also G\^{a}teaux-differentiable, and for $\psi\in C_0^{\infty}({\mathbb{R}}^N)$,
\begin{eqnarray}\label{def.J'}
\langle J'(v),\psi\rangle=\int_{{\mathbb{R}}^N}\nabla v\nabla\psi
-\lambda\int_{{\mathbb{R}}^N}c(x)h(v)h'(v)\psi-\int_{{\mathbb{R}}^N}f(h(v))h'(v)\psi.
\end{eqnarray}

Since $u=h(v)$, we have $\nabla u=h'(v)\nabla v$. Moreover, from (\ref{def.h}), we have $\nabla v=\beta(u)\nabla u$. We get immediately
that $h'(v)=[\beta(u)]^{-1}=[\beta(h(v))]^{-1}$ and $\nabla u=[\beta(h(v))]^{-1}\nabla v$.
Now assume that $v\in X$, $v>0$ such that equality $\langle J'(v),\psi\rangle=0$ holds. We choose $\psi=\beta(h(v))\varphi$, then
$\nabla \psi=\beta(h(v))\nabla\varphi+\beta'(h(v))h'(v)\varphi\nabla v$. Let $u=h(v)$, then $u\in X$. Moreover, we have $h'(v)\psi=\varphi$ and
\begin{eqnarray*}
\nabla v\nabla\psi&=&\beta(h(v))\nabla v\nabla\varphi+\beta'(h(v))h'(v)\varphi|\nabla v|^2 \\
    &=&\nabla u\nabla \varphi+\beta(u)\beta'(u)\varphi|\nabla u|^2.
\end{eqnarray*}
Thus from (\ref{def.J'}), we obtain that
\[
\int_{{\mathbb{R}}^N}\beta^2(u)\nabla u\nabla\varphi+\int_{{\mathbb{R}}^N}\beta(u)\beta'(u)\varphi|\nabla u|^2
-\lambda\int_{{\mathbb{R}}^N}c(x)u\varphi-\int_{{\mathbb{R}}^N}f(u)\varphi=0.
\]
This implies that $u$ such that (\ref{def.I'}) holds.
In summary, in order to find a weak solution to Eq.(\ref{eq.sch5-3}), it suffices to find a weak solution to the following equaiton
\begin{eqnarray}\label{eq.sch5-h}
-\Delta v-\lambda c(x)h(v)h'(v)=f(h(v))h'(v),\quad x\in{\mathbb{R}}^N.
\end{eqnarray}

We denote
\[
      \tilde{r}=\left\{
             \begin{array}{ll}
               +\infty, & \hbox{$0<\alpha\leq\frac{1}{2^*}$;} \\
               1+\frac{1}{2^*\alpha-1}, & \hbox{$\frac{1}{2^*}<\alpha<\frac{1}{2}$.}
             \end{array}
           \right.
\]
We assume that
\begin{description}
  \item[(c)] the function $c(x)\in C({\mathbb{R}^N},{\mathbb{R}})$, $c(x)>0$ and there exists $r\in(\frac{N}{2},\tilde{r})$ such that $c(x)\in L^{r}(\mathbb{R}^N)$.
  \item[(f)] assume that $q\in(2,2^*)$ and either \\
(i) $\frac{1}{4}<\alpha<\frac{1}{2}$, $q>\frac{4}{N-2}+4\alpha$ or \\
(ii) $0<\alpha\leq\frac{1}{4}$, $q>\frac{N+2}{N-2}$ holds.
\end{description}

Let $\lambda^*=2\alpha S\|c\|_{r}^{-1}$, where $S$ is the best constant for the Sobolev imbedding $H^1(\mathbb{R}^N)$ into $L^{2^*}(\mathbb{R}^N)$.
The following theorem is the main result of this paper.

\begin{thm}\label{thm.m1}
Assume that (c) and (f) hold,
then for $\lambda\in[0,\lambda^*)$, problem (\ref{eq.sch5-3}) has a positive weak solution $u\in X$.
Moreover, if $0<\alpha\leq {1}/{2^*}$, then $u\in H^1(\mathbb{R}^N)$.
\end{thm}

\begin{rem}\label{rem.m1-1}
We have $2^*-1<\frac{4}{N-2}+4\alpha<2^*$ for $\frac{1}{4}<\alpha<1/2$ and $\frac{N+2}{N-2}=2^*-1$.
\end{rem}

\begin{rem}\label{rem.m1-2}
Under the assumptions of Theorem \ref{thm.m1}, if $\|\lambda c\|_{r}< 2\alpha S$ for some $r\in(N/2,\tilde{r})$, then problem (\ref{eq.sch5-3}) has a positive solution in $X$.
\end{rem}

Let
\begin{eqnarray*}
&a_0=\frac{1+(2\alpha)^2}{1+2\alpha}\in(2\alpha,1),\quad\forall\alpha\in(0,1/2)\\
&a_1=\max\{a_0,2/2^*\},\quad \tilde{r}_1=2^*a_1/(2^*a_1-2).
\end{eqnarray*}

\begin{cor}\label{rem.m1-3}
Assume (f) holds and
\begin{description}
  \item[(c)$'$] the function $c(x)\in C({\mathbb{R}^N},{\mathbb{R}})$, $c(x)>0$ and there exists $r\in(N/2,\tilde{r}_1)$ such that $c(x)\in L^{r}(\mathbb{R}^N)$,
\end{description}
 then for $\lambda\in[0,S\|c\|_{r}^{-1})$, problem (\ref{eq.sch5-3}) has a positive solution.
\end{cor}

\begin{rem}\label{rem.m1-4}
By H\"{o}lder's inequality and Sobolev's inequality, we can prove that $\lambda^*=S\|c\|_{r}^{-1}\leq\lambda_1$, where $\lambda_1$ is the first eigenvalue of the following equation
\[
-\Delta u=\lambda c(x)u,\quad x\in{\mathbb{R}^N},
\]
that is,
\[
\lambda_1=\inf_{u\in X\setminus\{0\}}\frac{\int_{\mathbb{R}^N}|\nabla u|^{2}\mbox{d}x}{\int_{\mathbb{R}^N}c(x)u^2\mbox{d}x}.
\]
It is worthy of pointing out that computing the value $\lambda^*$ is much easier than obtaining $\lambda_1$.
Moreoer, the assumptions allow $c(x)$ to belong to a wide class of function space.
\end{rem}

The proof of Theorem \ref{thm.m1} is also applicable to problems at subcritical growth.

Let us consider the following equation
\begin{eqnarray}\label{eq.sch5-4}
-\Delta u-\lambda c(x)u-\kappa\alpha(\Delta(|u|^{2\alpha}))|u|^{2\alpha-2}u
    =|u|^{q-2}u,\quad x\in{\mathbb{R}}^N.
\end{eqnarray}

\begin{thm}\label{thm.m3}
Assume that (c) holds,
then for $\lambda\in[0,\lambda^*)$, $q_0<q<2^*$, $q_0=\max\{2,2^*(2\alpha)\}$, problem (\ref{eq.sch5-4}) has a positive weak solution $u\in X$.
Moreover, if $0<\alpha\leq {1}/{2^*}$, then $u\in H^1(\mathbb{R}^N)$.
\end{thm}

In Sect.\ref{sec:5-2}, we first study the properties of the function $h$; then we prove that the functional $J$ is well defined on $X$, continuous in $X$ and G\^{a}teaux-differentiable in $X$, see Proposition \ref{prop.J}. These are crucial steps since we only have $D^{1,2}(\mathbb{R}^N)\hookrightarrow L^{2^*}(\mathbb{R}^N)$. To this end, we establish several imbedding results. In the end of the section, we show that $J$ has the mountain pass geometry.
In Sect.\ref{sec:5-3}, we prove that every Palais-Smale ((PS) in short) sequence $\{v_n\}$ of $J$ is bounded in $X$ and analyze the properties of (PS) sequence. In Sect.\ref{sec:5-4}, we employ the mountain pass theorem without Palais-Smale condition to prove the existence of positive solution to Eq.(\ref{eq.sch5-h}). This involving the computation of a mountain pass level $c\in(0,\frac{1}{N}S^{N/2})$.

In this paper, $\|\cdot\|_p$ denotes the norm of $L^p(\mathbb{R}^N)$; $C,~C_1,~C_2,\cdots$ denote positive constants.


\section{ Mountain pass geometry}
\setcounter{equation}{0}
\label{sec:5-2}

In this section, we will give some properties of the transformation $h$
and establish the geometric hypotheses of the mountain pass geometry.

\begin{lem}\label{lem.h}
The function $h(t)$ has the following properties,
\begin{description}
  \item[(1)] $h(t)$ is odd, invertible, increasing and of class $C^1$ for $0<\alpha<1/2$, of class $C^2$ for $0<\alpha\leq1/4$;
  \item[(2)] $|h'(t)|\leq1$ for all $t\in{\mathbb{R}}$;
  \item[(3)] $|h(t)|\leq |t|$ for all $t\in{\mathbb{R}}$;
  \item[(4)] $h^{2\alpha}(t)/t\to\sqrt{2/\kappa}$ as $t\to0^+$;
  \item[(5)] $2\alpha h(t)\leq 2\alpha h'(t)t\leq h(t)$ for $t>0$;
  \item[(6)] $h(t)/t\to1$ as $t\to+\infty$;
  \item[(7)] $|h(t)|$ is convex;
  \item[(8)] $h'(t)\leq (2\kappa\alpha^2)^{-\vartheta/2}|h(t)|^{(1-2\alpha)\vartheta}$ for all $\vartheta\in[0,1].$
\end{description}
\end{lem}
{\bf Proof}\quad
For part (1), $h(t)$ is odd and invertible by definition. Since $h'(t)=[\beta(h(t))]^{-1}\in(0,1)$, we have $h(t)$ is increasing and of class $C^1$
for $0<\alpha<1/2$. By direct computation, we have
\[
h''(t)=2\kappa\alpha^2(1-2\alpha)\frac{|h(t)|^{-4\alpha}h(t)}{\big(2\kappa\alpha^2+|h(t)|^{2(1-2\alpha)}\big)^2},
\]
so we have $h(t)$ is of class $C^2$ for $0<\alpha\leq1/4$.

Part (2) is obvious. For part (3), assume that $t>0$ and note that $\beta(h(t))>1$, we have
\begin{eqnarray*}\label{lem.h-1}
h(t)=\int_0^{h(t)}ds\leq \int_0^{h(t)}\beta(s)ds=t.
\end{eqnarray*}
Then part (3) follows since $h$ is odd.

For part (4), note that from part (3) we have $h(t)\to 0$ as $t\to 0$. Thus we can employ
L'Hospital's principle to prove that
\[
\lim_{t\to0^+}\frac{h^{2\alpha}(t)}{t}=\lim_{t\to0^+}2\alpha h^{2\alpha-1}(t)h'(t)=\sqrt{\frac{2}{\kappa}}.
\]

For part (5), we prove the right-hand side inequality. Let $H(t)=h(t)/h'(t)$ and $\tilde{H}(t)=H(t)-2\alpha t$. Then $\tilde{H}(0)=0$.
We prove that $\tilde{H}'(t)\geq0$, i.e. $H'(t)\geq 2\alpha$, and this implies the conclusion.
In fact, note that $h(t)$ has same sign of $t$, for $t=0$, by part (4) we have
\[
H'(0)=\lim_{t\to0}\frac{H(t)}{t}=\lim_{t\to0}\sqrt{\frac{2}{\kappa}}\frac{|H(t)|}{|h(t)|^{2\alpha}}=\sqrt{\frac{2}{\kappa}}\sqrt{2\kappa\alpha^2}=2\alpha.
\]
For $t\neq0$, we have
\begin{eqnarray*}
H'(t)&=&\bigg(\frac{h(t)\big(2\kappa\alpha^2+|h(t)|^{2(1-2\alpha)}\big)^{1/2}}{|h(t)|^{1-2\alpha}}\bigg)' \\
&\geq&\frac{|h(t)|^{2(1-2\alpha)}-(1-2\alpha)h(t)(2\kappa\alpha^2+|h(t)|^{2(1-2\alpha)})^{1/2}|h(t)|^{-1-2\alpha}h(t)h'(t)}{|h(t)|^{2(1-2\alpha)}} \\
&=&\frac{|h(t)|^{2(1-2\alpha)}-(1-2\alpha)|h(t)|^{2(1-2\alpha)}}{|h(t)|^{2(1-2\alpha)}} \\
&=&2\alpha.
\end{eqnarray*}
The left-hand side inequality can be proved similarly.

For part (6), we have $h'(t)>1/2$ for $t>0$ sufficiently large. So we conclude that $h(t)\to+\infty$ as $t\to+\infty$.
Thus by employing Hospital principle again, we have $\lim_{t\to+\infty}h(t)/t=\lim_{t\to+\infty}h'(t)=1$.

(7) Note that $h(t)$ is odd and $h(0)=0$, $h''(t)>0$ for $t>0$, we conclude that $|h(t)|$ is convex.

(8) For any $\vartheta\in[0,1]$, by the definition of $h'(t)$ and note that $0\leq h'(t)\leq1$, we have
\begin{eqnarray*}
h'(t)=h'(t)^{\vartheta}h'(t)^{1-\vartheta}
    \leq h'(t)^{\vartheta}
    \leq (2\kappa\alpha^2)^{-\vartheta/2}|h(t)|^{(1-2\alpha)\vartheta}.
\end{eqnarray*}
This ends the proof of the lemma.
$\quad\Box$

\begin{lem}\label{lem.h'-1}
For $t,~s\in{\mathbb{R}}$, we have
\[
|h'(t)-h'(s)|\leq 2^{1-\vartheta}(2\kappa\alpha^2)^{-\vartheta/2}|h(t)-h(s)|^{(1-2\alpha)\vartheta},\quad\forall~\vartheta\in[0,1].
\]
\end{lem}
{\bf Proof}\quad Consider $\eta(t):=t^{a},~t,~s>0,~a\in(0,1)$. We have
\begin{eqnarray}\label{lem.h'-1-1}
\eta(|t-s|)\geq|\eta(t)-\eta(s)|.
\end{eqnarray}
Now for $t,~s\in{\mathbb{R}}$, by direct computation and using (\ref{lem.h'-1-1}), we have
\begin{eqnarray*}\label{lem.h'-1-2}
H(t,s)&:=&|h'(t)-h'(s)| \\
&\leq&\frac{({2\kappa\alpha^2})^{1/2}||h(t)|^{1-2\alpha}-|h(s)|^{1-2\alpha}|}{(2\kappa\alpha^2+|h(t)|^{2(1-2\alpha)})^{1/2}(2\kappa\alpha^2+|h(s)|^{2(1-2\alpha)})^{1/2}}\\
&\leq& ({2\kappa\alpha^2})^{-1/2}|h(t)-h(s)|^{1-2\alpha}.
\end{eqnarray*}
Note that by (2) of Lemma \ref{lem.h}, we also have $H(t,s)\leq2$, we conclude that
\[
H(t,s)=H^{1-\vartheta}(t,s)H^{\vartheta}(t,s)\leq2^{1-\vartheta}(2\kappa\alpha^2)^{-\vartheta/2}|h(t)-h(s)|^{(1-2\alpha)\vartheta}
\]
for all $\vartheta\in[0,1]$.
$\quad\Box$

\begin{lem}\label{lem.embed-2}
If $v\in X$, then $h(v)\in L^r(\mathbb{R}^N),~r\in[2^*(2\alpha),2^*]$.
\end{lem}
{\bf Proof}\quad Let $v\in X$, then we have $\|\nabla v\|_2\leq C$ for some constant $C>0$.
For $a\in[2\alpha,1]$, we have
\begin{eqnarray}\label{lem.embed-2-0}
\nabla(|h(v)|^a)=a|h(v)|^{a-2}h(v)h'(v)\nabla v=\frac{a|h(v)|^{a-2}h(v)|h(v)|^{1-2\alpha}}{(2\alpha+|h(v)|^{2(1-2\alpha)})^{1/2}}\nabla v
\end{eqnarray}
In order to ensure that $\|\nabla(h(v)^a)\|_2\leq C$, it requires that
\[
0\leq (a-1)+(1-2\alpha)=a-2\alpha\leq 1-2\alpha,
\]
that is, $2\alpha\leq a\leq 1$. Since
\begin{eqnarray}\label{lem.embed-2-1}
|\nabla(|h(v)|^a)|^2\leq\frac{1}{2\alpha}|\nabla v|^2,\quad \forall a\in[2\alpha,1],
\end{eqnarray}
we obtain that $\|\nabla(|h(v)|^a)\|_2\leq C/(2\alpha)$ for all $a\in[2\alpha,1]$. Similar to the proof of Sobolev's inequality, we have
$|h(v)|^a\in L^{2^*}(\mathbb{R}^N)$, this gives that $h(v)\in L^{2^*a}(\mathbb{R}^N)$ for $a\in[2\alpha,1]$.
$\quad\Box$

\begin{lem}\label{lem.embed-3}
The map: $v\to h(v)$ from $X$ into $L^r({\mathbb{R}}^N)$ is continuous for $2^*(2\alpha)\leq r\leq 2^*$.
\end{lem}
{\bf Proof}\quad
For any $r\in[2^*(2\alpha),2^*]$, there exists $a\in[2\alpha,1]$ such that $r=2^*a$. Then for any sequence $\{v_n\}\subset X$ that converges strongly to $0$, by Sobolev's inequality and (\ref{lem.embed-2-1}), we have
\begin{eqnarray}\label{lem.embed-3-1}
\|h(v_n)\|_{r}^{2a}=\||h(v_n)|^a\|_{2^*}^{2}\leq S^{-1}\|\nabla |h(v_n)|^a\|_{2}^{2}\leq (2\alpha S)^{-1}\|v_n\|_X^2,
\end{eqnarray}
which tends to $0$ as $n\to\infty$.
Thus the map $v\to h(v)$ from $X$ into $L^{r}({\mathbb{R}}^N)$ is continuous.
$\quad\Box$

By using (2)-(3) of Lemma \ref{lem.h}, we can also prove that

\begin{lem}\label{lem.embed-1}
The map: $v\to h(v)$ from $H^1({\mathbb{R}}^N)$ into $L^p({\mathbb{R}}^N)$ is continuous for $2\leq p\leq 2^*$,
and is locally compact for $2\leq p<2^*$.
\end{lem}
{\bf Proof}\quad
The proof is similar to that for Lemma \ref{lem.embed-3}.
Firstly, assume that $\{v_n\}\subset H^1({\mathbb{R}}^N)$ converges strongly to $0$, then by part (3) of Lemma \ref{lem.h} and Sobolev's inequality, we have
\begin{eqnarray*}\label{lem.embed-1-1}
\|h(v_n)\|^2_{2^*}\leq \|v_n\|^2_{2^*}\leq S^{-1}\|v_n\|^2_X,
\end{eqnarray*}
which tends to $0$ as $n\to\infty$.
Thus the map $v\to h(v)$ from $H^1({\mathbb{R}}^N)$ into $L^{2^*}({\mathbb{R}}^N)$ is continuous. By interpolation inequality, we obtain that
the map $v\to h(v)$ from $H^1({\mathbb{R}}^N)$ into $L^{p}({\mathbb{R}}^N)$ is continuous for $2\leq p\leq2^*$.

Secondly, assume that $\{v_n\}$ is bounded in $H^1({\mathbb{R}}^N)$, exist a subsequence (we still denote it by $\{v_n\}$) and a $v\in H^1({\mathbb{R}}^N)$ such that
$v_n\to v$ locally in $L^{p}({\mathbb{R}}^N)$ for $2\leq p<2^*$.
By mean value theorem and (2)-(3) of Lemma \ref{lem.h}, we have
\begin{eqnarray}\label{lem.embed-1-2}
|h(v_n)-h(v)|\leq|h'(v+\theta(v_n-v))||v_n-v|\leq |v_n-v|,
\end{eqnarray}
where $\theta\in(0,1)$. Thus the map $v\to h(v)$ from $H^1({\mathbb{R}}^N)$ into $L^{p}({\mathbb{R}}^N)$ is locally compact for $2\leq p<2^*$.
$\quad\Box$

Noting  that $H^1({\mathbb{R}}^N)\subset D^{1,2}({\mathbb{R}}^N)$ is dense, we thus can combine Lemma \ref{lem.embed-3} and Lemma \ref{lem.embed-1} to conclude that
\begin{cor}\label{lem.embed-4}
The map: $v\to h(v)$ from $H^1({\mathbb{R}}^N)$ into $L^p({\mathbb{R}}^N)$ is continuous for $2_{*}\leq p\leq 2^*$,
and is locally compact for $2_{*}\leq p<2^*$, where $2_{*}=\min\{2,2^*(2\alpha)\}$.
\end{cor}
{\bf Proof}\quad
If $2\leq 2^*(2\alpha)$, then the conclusions hold obviously. Now we assume that
$2^*(2\alpha)<2$ and $p\in[2^*(2\alpha),2)$. Then the first part of the corollary follows from (\ref{lem.embed-3-1}) and the fact that $\|v_n\|_X^2\leq\|v_n\|_{H^1(\mathbb{R}^N)}^2$. On the other hand, by H\"{o}lder's inequality, we have
\begin{eqnarray}\label{lem.embed-4-1}
\|h(v_n)-h(v)\|_{p}\leq\|h(v_n)-h(v)\|_{2^*(2\alpha)}^{\theta}\|h(v_n)-h(v)\|_{2}^{1-\theta}
\end{eqnarray}
for some $\theta\in[0,1]$. Then the second part of the corollary holds by Lemma \ref{lem.embed-1}.
$\quad\Box$

\begin{lem}\label{lem.embed-5}
The map: $v\to h(v)$ from $X$ into $L^r({\mathbb{R}}^N)$ is locally compact for $2^*(2\alpha)\leq r<2^*$.
\end{lem}
{\bf Proof}\quad
We should only prove that the conclusion holds for the case $2^*(2\alpha)<2$ and $r\in[2^*(2\alpha),2)$. If $2^*(2\alpha)\leq2$, by Lemma \ref{lem.embed-2}, for any $v\in X$, we have $h(v)\in L^2(\mathbb{R}^N)$, thus $h(v)\in H^1(\mathbb{R}^N)$.
Now assume that $\{v_n\}\subset X$ is bounded, by (\ref{lem.embed-2-1}), $\{h(v_n)\}\subset L^2(\mathbb{R}^N)$ is bounded uniformly in $n$.
Thus $\{h(v_n)\}$ is bounded in $H^1(\mathbb{R}^N)$. Then by Lemma \ref{lem.embed-1}, we conclude that $\{h(v_n)\}$ is locally compact in $L^2(\mathbb{R}^N)$. Finally, we obtain that the map $v\to h(v)$ from $X$ into $L^{r}({\mathbb{R}}^N)$ is locally compact for $r\in[2^*(2\alpha),2)$ by inequality (\ref{lem.embed-4-1}).
$\quad\Box$

If $0<\alpha\leq {1}/{2^*}$, then we have $2^*(2\alpha)\leq2$. Thus, according to the proof of Lemma \ref{lem.embed-5}, we have

\begin{thm}\label{thm.u-in-H1}
Assume that $0<\alpha\leq {1}/{2^*}$, then every solution for (\ref{eq.sch5-3}) belongs to $H^1(\mathbb{R}^N)$.
\end{thm}

Let $L^2(\mathbb{R}^N,c(x))$ be the weighted Lebesgue space defined by
\[
L^2(\mathbb{R}^N,c(x)):=\bigg\{u\in L^2(\mathbb{R}^N):\int_{\mathbb{R}^N}c(x)u^2\mbox{d}x<+\infty\bigg\},
\]
and endowed with the norm
\[
\|u\|_{2,c}:=\bigg(\int_{\mathbb{R}^N}c(x)u^2\mbox{d}x\bigg)^{1/2}.
\]

We have
\begin{lem}\label{lem.embed-Lc}
Assume hypothesis (c) holds, then the map $v\rightarrow h(v)$ from $X$ into $L^2(\mathbb{R}^N,c(x))$ is continuous and compact.
\end{lem}
{\bf Proof}\quad
Let $s>0$ satisfies that $\frac{2}{s}+\frac{1}{r}=1$, where $r\in(\frac{N}{2},\tilde{r})$ and $\tilde{r}$ is given by assumption (c).
By Lemma \ref{lem.embed-2}, we have $c(x)\in L^{r}({\mathbb{R}}^N)$ and $h(v)\in L^{s}({\mathbb{R}}^N)$.
Then for $v\in X$, by H\"{o}lder's inequality, we have
\begin{eqnarray}\label{lem.embed-Lc-1}
\int_{{\mathbb{R}}^N}c(x)h(v)^2\mbox{d}x\leq \|c\|_{r'}\|h(v)\|_{s}^2<+\infty.
\end{eqnarray}
By Lemma \ref{lem.embed-3}, we obtain that the the map $v\to h(v)$ from $X$ into $L^2(\mathbb{R}^N,c(x))$ is continuous.

Now assume that $\{v_n\}\subset X$ is bounded, then there exist a subsequence (still denoted by $\{v_n\}$), and a $v\in X$, such that $v_n\rightharpoonup v$ in $X$, $v_n\to v$ in $L_{\rm loc}^t({\mathbb{R}}^N),~1\leq t<2^*$.
On the other hand, since by (\ref{lem.embed-2-1}), $\{h(v_n)\}\subset X$ is also bounded, there exist a subsequence (still denoted by $\{h(v_n)\}$), and a $w\in X$, such that $h(v_n)\rightharpoonup w$ in $X$, $h(v_n)\to w$ in $L_{\rm loc}^t({\mathbb{R}}^N),~2^*(2\alpha)\leq t<2^*$.

We claim that $w=h(v)$ a.e. in ${\mathbb{R}}^N$. Indeed, for any $\delta>0$, there exists $R_\delta>0$ such that
$\int_{B_{R_\delta}^c}c(x)h(v)^2<\delta/3$ and $\int_{B_{R_\delta}^c}c(x)w^2<\delta/3$, where $B_{R_\delta}=\{x\in{\mathbb{R}}^N:|x|\leq R_\delta\}$ and $B_{R_\delta}^c={\mathbb{R}}^N\setminus B_{R_\delta}$. Thus, by (\ref{lem.embed-1-2}) and the locally compact imbedding, we have for $n$ sufficiently large,
\begin{eqnarray*}
0&\leq&\int_{{\mathbb{R}}^N}c(x)|h(v)-w|^2 \\
&=&\int_{B_{R_\delta}}c(x)|h(v)-w|^2+\int_{B_{R_\delta}^c}c(x)|h(v)-w|^2 \\
&\leq&\int_{B_{R_\delta}}c(x)|h(v_n)-h(v)|^2+\int_{B_{R_\delta}}c(x)|h(v_n)-w|^2 \\
&&+\int_{B_{R_\delta}^c}c(x)h(v)^2+\int_{B_{R_\delta}^c}c(x)w^2 \\
&<&\delta.
\end{eqnarray*}
This proves the claim.

Now for any $\varepsilon>0$, there exists $R_{\varepsilon}>0$ such that
\[
\int_{B_{R_\varepsilon}^c}|c(x)|^{r}\mbox{d}x\leq\bigg(\frac{\varepsilon}{4C}\bigg)^{r}
\]
where $C>0$ satisfies that $(2\alpha S)^{-1}\|v_n\|_{X}^2\leq C$. Thus by H\"{o}lder's inequality and (\ref{lem.embed-2-1}), we have
\[
\int_{B_{R_\varepsilon}^c}c(x)|h(v_n)-h(v)|^{2}\mbox{d}x
\leq\|c\|_{L^{r}(B_{R_\varepsilon}^c)}\|h(v_n)-h(v)\|^2_{L^{s}(B_{R_\varepsilon}^c)}\leq\frac{\varepsilon}{2}
\]
Since $X\hookrightarrow L^2(B_R)$ is compact, it follows that
there exists $n_0\in{\mathbb{N}}$ such that for $n\geq n_0$,
\[
\int_{B_{R_\varepsilon}}c(x)|h(v_n)-h(v)|^2\mbox{d}x\leq\frac{\varepsilon}{2}.
\]
Thus we obtain that
\[
\int_{{\mathbb{R}}^N}c(x)|h(v_n)-h(v)|^2\mbox{d}x\leq\bigg(\int_{B_{R_\varepsilon}}+\int_{B_{R_\varepsilon}^c}\bigg)c(x)|h(v_n)-h(v)|^2\mbox{d}x\leq{\varepsilon}.
\]
This implies that the map $v\rightarrow h(v)$ from $X$ into $L^2(\mathbb{R}^N,c(x))$ is compact.
$\quad\Box$

Now we come back to the discussion of functional $J$. We have

\begin{prop}\label{prop.J}
Under assumptions of Theorem \ref{thm.m1}, the functional $J$ has the following properties:$\\$
(1) $J$ is well defined on $X$. $\\$
(2) $J$ is continuous in $X$. $\\$
(3) $J$ is G\^{a}teaux-differentiable.
\end{prop}
{\bf Proof}\quad
(1) Firstly, for $v\in X$, by Lemma \ref{lem.embed-Lc}, we have $\int_{\mathbb{R}^N}c(x)h(v)^2\mbox{d}x<+\infty$.
Next, by assumptions of Theorem \ref{thm.m1}, for $\alpha\in(\frac{1}{4},\frac{1}{2})$, we have $q>\frac{4}{N-2}+4\alpha>2^*(2\alpha)$;
for $\alpha\in(0,\frac{1}{4}]$, we have $q>2^*-1> 2^*(2\alpha)$. Then Lemma \ref{lem.embed-3} implies that $\|h(v)\|_{q}^{q}<+\infty$.
Thus we have
\[
\int_{{\mathbb{R}}^N}F(x,h(v))\mbox{d}x\leq C(\|h(v)\|_{q}^{q}+\|h(v)\|_{2^*}^{2^*})<+\infty.
\]
These show that $J$ is well defined on $X$.

(2) Assume that $v_n\to v$ in $X$. By Sobolev's inequality, $\|v_n-v\|_{2^*}\to0$. By mean value theorem, (8) of Lemma \ref{lem.h}  and H\"{o}lder's inequality, we have
\begin{eqnarray*}
&&\bigg|\int_{{\mathbb{R}}^N}c(x)[h(v_n)^2-h(v)^2]\mbox{d}x\bigg| \\
&&\quad =\bigg|\int_{{\mathbb{R}}^N}2c(x)h(v+\theta_n(v_n-v))h'(v+\theta_n(v_n-v))(v_n-v)\mbox{d}x\bigg| \\
&&\quad \leq 2(2\alpha)^{-\vartheta_r/2}\int_{{\mathbb{R}}^N}c(x)|h(v+\theta_n(v_n-v))|^{1+(1-2\alpha)\vartheta_r}|v_n-v|\mbox{d}x \\
&&\quad \leq 2(2\alpha)^{-\vartheta_r/2}\|c\|_r\|h(v+\theta_n(v_n-v))\|_{s}^{1+(1-2\alpha)\vartheta_r}\|v_n-v\|_{2^*}
\to0,
\end{eqnarray*}
where $\theta_n\in(0,1)$ and $s=[1+(1-2\alpha)\vartheta_r](1-\frac{1}{2^*}-\frac{1}{r})^{-1}$ with some $\vartheta_r\in[0,1]$ such that $s\in[2^*(2\alpha),2^*]$. Here, in the last inequality, we have used Lemma \ref{lem.embed-3} to obtain that $\|h(v+\theta_n(v_n-v))\|_{s}<+\infty$. Likewise, together with (2)-(3) of Lemma \ref{lem.h}, for $q\geq2^*(2\alpha)$,
\begin{eqnarray*}
&&\bigg|\int_{{\mathbb{R}}^N}[F(h(v_n))-F(h(v))]\mbox{d}x\bigg| \\
&&\quad \leq\int_{{\mathbb{R}}^N}\big|f(h(v+\theta_n(v_n-v)))h'(v+\theta_n(v_n-v))(v_n-v)\big|\mbox{d}x \\
&&\quad \leq\int_{{\mathbb{R}}^N}\big|h(v+\theta_n(v_n-v))\big|^{q-1}h'(v+\theta_n(v_n-v))\big|(v_n-v)\big|\mbox{d}x \\
&&\qquad    +\int_{{\mathbb{R}}^N}\big|v+\theta_n(v_n-v)\big|^{2^*-1}\big|(v_n-v)\big|\mbox{d}x \\
&&\quad \leq (2\alpha)^{-\vartheta_q/2}\big(\|h(v+\theta_n(v_n-v))\|_{t}^{(q-1)+(1-2\alpha)\vartheta_q}\|v_n-v\|_{2^*}\big) \\
&&\qquad    +2^{2^*-2}\big(\|v\|_{2^*}^{2^*-1}\|v_n-v\|_{2^*}+\|v_n-v\|_{2^*}^{2^*}\big)
\to0,
\end{eqnarray*}
where $t=[(q-1)+(1-2\alpha)\vartheta_q]2^*(2^*-1)^{-1}$ with some $\vartheta_q\in[0,1]$ such that $t\in[2^*(2\alpha),2^*]$.
Thus we have $J(v_n)\to J(v)$, that is, $J$ is continuous in $X$.

(3) Since $h\in C^1(\mathbb{R})$, for $v\in X$, $t>0$ and for any $\phi\in X$, by mean value theorem, we have
\[
\frac{1}{t}\int_{{\mathbb{R}}^N}c(x)\big[h(v+t\phi)^2-h(v)^2\big]\mbox{d}x=\int_{{\mathbb{R}}^N}2c(x)h(v+\theta t\phi)h'(v+\theta t\phi)\phi\mbox{d}x,
\]
where $\theta\in(0,1)$. By mean value theorem, we have
\begin{eqnarray*}
&& I:=\bigg|\int_{{\mathbb{R}}^N}c(x)h(v+\theta t\phi)h'(v+\theta t\phi)\phi\mbox{d}x-\int_{{\mathbb{R}}^N}c(x)h(v)h'(v)\phi\mbox{d}x\bigg| \\
&&\quad \leq\int_{{\mathbb{R}}^N}c(x)\big|h(v+\theta t\phi)-h(v)\big||h'(v+\theta t\phi)||\phi|\mbox{d}x \\
&&\qquad            +\int_{{\mathbb{R}}^N}c(x)\big|h(v)\big|\big|h'(v+\theta t\phi)-h'(v)\big||\phi|\mbox{d}x \\
&&\quad \leq \theta t\int_{{\mathbb{R}}^N}c(x)\big|h'(v+\xi\theta t\phi)\big||h'(v+\theta t\phi)||\phi|^2\mbox{d}x \\
&&\qquad            +\int_{{\mathbb{R}}^N}c(x)\big|h(v)\big|\big|h'(v+\theta t\phi)-h'(v)\big||\phi|\mbox{d}x
:=I_1+I_2,
\end{eqnarray*}
where $\theta,~\xi\in(0,1)$.

(i) We consider $I_1$.
If $r=\frac{N}{2}$, then by (2) of Lemma \ref{lem.h} and H\"{o}lder's inequality,
\[
I_1\leq \theta t\int_{{\mathbb{R}}^N}c(x)|\phi|^2\mbox{d}x \leq \theta t\|c\|_{N/2}\|\phi\|_{2^*}^2\to0,\quad t\to0.
\]
Otherwise, let us consider $s_1:=s_1(r)=(1-\frac{2}{2^*}-\frac{1}{r})^{-1}$ defined in $r\in(\frac{N}{2},\tilde{r})$. We have $s_1(r)$ is decreasing and $s_1\in(\frac{2^*\alpha}{1-2\alpha},+\infty)$ for $\alpha\in(\frac{1}{4},\frac{1}{2})$, $s_1\in(\frac{2^*}{2^*-2},+\infty)$ for $\alpha\in(0,\frac{1}{4}]$. Let $t_1:=t_1(s_1,\vartheta_1)=2s_1(1-2\alpha)\vartheta_{1}$. Note that for any ${r}\in(\frac{N}{2},\tilde{r})$, there exists $\vartheta_{1}:=\vartheta_1(r)\in[0,1]$ such that $t_1=2^*(2\alpha)$,
by the definition of $h'(t)$ and H\"{o}lder's inequality,
\begin{eqnarray*}
I_1\leq \theta t(2\alpha)^{-\vartheta_{1}}\|c\|_{r}\|h(v+\xi\theta t\phi)\|_{2^*(2\alpha)}^{(1-2\alpha)\vartheta_{1}}\|h(v+\theta t\phi)\|_{2^*(2\alpha)}^{(1-2\alpha)\vartheta_{1}}\|\phi\|_{2^*}^2
\to0,
\quad t\to0.
\end{eqnarray*}

(ii) We consider $I_2$. Firstly, for $\alpha\in(\frac{1}{2^*},\frac{1}{2})$, we have $2^*(2\alpha)\in(2,2^*)$. Let $s_2:=s_2(r)=(1-\frac{1}{2^*}-\frac{1}{r})^{-1}$,
then for $r\in[\frac{N}{2},\frac{2^*(2\alpha)}{(2^*-1)2\alpha-1}]$, we have $s_2\in[2^*(2\alpha),2^*]$.
By H\"{o}lder's inequality and Lebesgue's dominated convergence theorem,
\[
I_2\leq \|c\|_{r}\|h(v)\|_{s_2}\|(h'(v+\theta t\phi)-h'(v))\phi\|_{2^*}\to0,\quad t\to0;
\]
For $r\in(\frac{2^*(2\alpha)}{(2^*-1)2\alpha-1},\frac{2^*\alpha}{2^*\alpha-1})$, let $s_3:=s_3(r)=(1-\frac{1}{2^*}-\frac{1}{2^*(2\alpha)}-\frac{1}{r})^{-1}$. We have $s_3\in(\frac{2^*(2\alpha)}{1-2\alpha},+\infty)$.
Let $t_3:=t_3(s_3,\vartheta_3)=(1-2\alpha)\vartheta_3 s_3$, then for any $r\in(\frac{2^*(2\alpha)}{(2^*-1)2\alpha-1},\frac{2^*\alpha}{2^*\alpha-1})$,
there exists $\vartheta_3:=\vartheta_3(r)\in[0,1]$ such that $t_3=2^*(2\alpha)$. By Lemma \ref{lem.h'-1}, H\"{o}lder's inequality and Lebesgue's dominated convergence theorem,
\begin{eqnarray*}
I_2 \leq C\|c\|_{r}\|h(v)\|_{2^*(2\alpha)}\|h(v+\theta t\phi)-h(v)\|_{2^*(2\alpha)}^{(1-2\alpha)\vartheta_3}
 \|\phi\|_{2^*}\to0,\quad t\to0;
\end{eqnarray*}

Secondly, for $\alpha\in(0,\frac{1}{2^*}]$, we have $2^*(2\alpha)\in(0,2]$. Let $s_4:=s_4(r)=(1-\frac{1}{2^*}-\frac{1}{r})^{-1}$, then
for $r\in[\frac{N}{2},+\infty)$, we have $s_4\in(\frac{2^*}{2^*-1},2^*]\subset(1,2^*]$. If $s_4\geq2^*(2\alpha)$,
then
\[
I_2\leq \|c\|_{r}\|h(v)\|_{s_4}\|(h'(v+\theta t\phi)-h'(v))\phi\|_{2^*}\to0,\quad t\to0;
\]
If $s_4<2^*(2\alpha)$,
we let $t_4:=t_4(s_4, \theta_4)=(1-2\alpha)\theta_4\frac{2^*(2\alpha)s_4}{2^*(2\alpha)-s_4}$, then $t_4$ is increasing in $s_4$ and
$t_4\in(2^*(2\alpha)\theta_4,(+\infty)\theta_4)$.
Note that there exists $\vartheta_4:=\vartheta_4(r)\in[0,1]$ such that $t_4=2^*(2\alpha)$, we have
\begin{eqnarray*}
I_2&\leq& \|c\|_{r}\|h(v)(h'(v+\theta t\phi)-h'(v))\|_{s_4}\|\phi\|_{2^*} \\
&\leq& C\|c\|_{r}\|h(v)\|_{2^*(2\alpha)}\|h(v+\theta t\phi)-h(v)\|_{2^*(2\alpha)}^{(1-2\alpha)\vartheta_4}\|\phi\|_{2^*}
\to0,\quad t\to0.
\end{eqnarray*}

In summary, from (i)-(ii), we conclude that $I\to0$. This means that
\[
\frac{1}{t}\int_{{\mathbb{R}}^N}c(x)\big[h(v+t\phi)^2-h(v)^2\big]\mbox{d}x\to\int_{{\mathbb{R}}^N}2c(x)h(v)h'(v)\phi\mbox{d}x.
\]
Likewise, for $q\geq2^*(2\alpha)$, we have
\[
\frac{1}{t}\int_{{\mathbb{R}}^N}\big[F(h(v+t\phi))-F(h(v))\big]\mbox{d}x\to\int_{{\mathbb{R}}^N}f(h(v))h'(v)\phi\mbox{d}x.
\]
These imply that $J$ is G\^{a}teaux-differentiable. This completes the proof.
$\quad\Box$

In the following, we consider the existence of positive solutions of Eq.(\ref{eq.sch5-h}). From variational methods, we will study the positive critical points of the following functional
\[
J^+(v)=\frac{1}{2}\int_{{\mathbb{R}}^N}|\nabla v|^2\mbox{d}x
-\frac{\lambda}{2}\int_{{\mathbb{R}}^N}c(x)h(v)^2\mbox{d}x-\int_{{\mathbb{R}}^N}F(h(v)^+)\mbox{d}x.
\]
To avoid cumbersome notations, in the rest of this paper, we still denote $J^+(v)$ and $F(h(v)^+)$ by $J(v)$ and $F(h(v))$ respectively.

\begin{lem}\label{lem.mp-1}
There exist $\rho_0,~a_0>0$ such that $J(v)\geq a_0$ for all $\|v\|_{X}=\rho_0$.
\end{lem}
{\bf Proof}\quad
Let $s>0$ satisfies that $\frac{2}{s}+\frac{1}{r}=1$, where $r\in(N/2,\tilde{r})$ and $\tilde{r}$ is given by assumption (c).
Note that $|h(v)|\leq |v|$, by Sobolev's inequality, we have
\begin{eqnarray}\label{lem.mp-1-1}
J(v)&=& \frac{1}{2}\int_{\mathbb{R}^N}|\nabla v|^2\mbox{d}x
-\frac{\lambda}{2}\int_{\mathbb{R}^N}c(x)h(v)^2\mbox{d}x
-\int_{\mathbb{R}^N}F(h(v))\mbox{d}x \nonumber\\
&\geq& \frac{1}{2}\|v\|_{X}^2-\frac{\lambda}{2}\|c\|_{r}\|h(v)\|_{s}^{2}
-\frac{1}{q}\|h(v)\|_{q}^{q}
-\frac{1}{2^*}\|h(v)\|_{2^*}^{2^*} \nonumber\\
&\geq& \frac{1}{2}(1-\lambda C_1S^{-1}\|c\|_{r})\|v\|_{X}^2-C_2(\|v\|_{X}^q+\|v\|_{X}^{2^*}),
\end{eqnarray}
where $C_1=1/(2\alpha)$ according to (\ref{lem.embed-2-1}). Let $\lambda^*=C_1S^{-1}\|c\|_{r}$. Note that $2^*>q>2$, then for $\lambda\in(0,\lambda^*)$, there exist $\rho>0$ and $a_0>0$ such that
$J(v)\geq a_0$ for all $\|v\|_{X}=\rho$.
$\quad\Box$

\begin{lem}\label{lem.mp-2}
There exists $v\in X$ such that $J(v)<0$.
\end{lem}
{\bf Proof}\quad
Given $\varphi\in C_0^{\infty}(\mathbb{R}^N,[0,1])$ with $\mbox{supt}(\varphi)=\overline{B}_2$ and $\varphi(x)=1$ for $x\in B_1$.
Note that $\lim\limits_{t\to+\infty}h(t\varphi)/t\varphi=1$, we have $F(h(t\varphi))\geq\frac{1}{2}F(t\varphi)$ for $t\in\mathbb{R}$ large enough.
Then we have
\begin{eqnarray*}
J(t\varphi)\leq\frac{t^2}{2}\int_{\mathbb{R}^N}|\nabla\varphi|^2
-\frac{\lambda t^2}{4}\int_{B_1}c(x)|\varphi|^2
-\frac{t^q}{2q}\int_{B_1}|\varphi|^q
-\frac{t^{2^*}}{22^*}\int_{B_1}|\varphi|^{2^*}.
\end{eqnarray*}
Let $v=t_0\varphi$ with $t_0>0$ sufficiently large, we have $J(v)<0$.
$\quad \Box$


\section{ Analysis of (PS) conditions}
\setcounter{equation}{0}
\label{sec:5-3}

As a consequence of Lemma \ref{lem.mp-1} and Lemma \ref{lem.mp-2}, there exists a Palais-Smale sequence $\{v_n\}$ of $J$ at level $c$ with
\begin{eqnarray}\label{def.c}
c=\inf_{\gamma\in\Gamma}\sup_{t\in[0,1]}J(\gamma(t))>0,
\end{eqnarray}
where $$\Gamma=\{\gamma\in C([0,1],X):\gamma(0)=0,\gamma(1)\neq0,J(\gamma(1))<0\}.$$
That is, $\{v_n\}$ satisfies $J(v_n)\to c,~J'(v_n)\to0$ as $n\to\infty$.

\begin{prop}\label{prop.ps-1}
Every Palais-Smale sequence $\{v_n\}$ for $J$ is bounded in $X$.
\end{prop}
{\bf Proof}\quad
Since $\{v_n\}\subset X$ is a Palais-Smale sequence, we have
\begin{eqnarray}\label{prop.ps-1-1}
J(v_n)=\frac{1}{2}\int_{\mathbb{R}^N}|\nabla v_n|^2dx
-\frac{\lambda}{2}\int_{\mathbb{R}^N}c(x)h(v_n)^2dx
-\int_{\mathbb{R}^N}F(h(v_n))dx\to c,
\end{eqnarray}
and for any $\psi\in C_0^{\infty}(\mathbb{R}^N)$,
\begin{eqnarray*}\label{prop.ps-1-2}
J'(v_n)\psi&=&\int_{\mathbb{R}^N}\Big[\nabla v_n\nabla\psi
-\lambda c(x)h(v_n)h'(v_n)\psi
-f(h(v_n))h'(v_n)\psi\Big]dx\nonumber\\
&=&o(1)\|\psi\|_X.
\end{eqnarray*}
Note that $h(t)/h'(t)\to0$ as $t\to0$, we have $h(t)/h'(t)\in X$ by direct computation. Moreover, since $C_0^{\infty}(\mathbb{R}^N)$ is dense in $X$, we can take $\psi=h(v_n)/h'(v_n)$ as test functions and get
\begin{eqnarray}\label{prop.ps-1-3}
\langle J'(v_n),\psi\rangle &=& \int_{{\mathbb{R}}^N}|\nabla v_n|^2-\lambda\int_{{\mathbb{R}}^N}c(x)h(v_n)^2-\int_{{\mathbb{R}}^N}f(h(v_n))h(v_n) \nonumber\\
&& -\int_{{\mathbb{R}}^N}\frac{2\alpha(1-2\alpha)}{2\alpha+|h(v_n)|^{2(1-2\alpha)}}|\nabla v_n|^2.
\end{eqnarray}
It follows that
\begin{eqnarray*}
c+o(1)\|v_n\|_{X} &=& J(v_n)-\frac{1}{q}\langle J'(v_n),\psi\rangle \\
    &\geq& \bigg(\frac{1}{2}-\frac{1}{q}\bigg)\int_{{\mathbb{R}}^N}|\nabla v_n|^2
    -\lambda\bigg(\frac{1}{2}-\frac{1}{q}\bigg)\int_{{\mathbb{R}}^N}c(x)h(v_n)^2.
\end{eqnarray*}
Similar to (\ref{lem.mp-1-1}), we obtain $\{v_n\}$ is bounded in $X$.
Note that $|\nabla h(v_n)|\leq|\nabla v_n|$, we conclude that $\{h(v_n)\}$ is also bounded in $X$.
$\quad\Box$

Since $v_n$ is a bounded Palais-Smale sequence, there exists $v\in X$ such that $v_n\rightharpoonup v$ in $X$. We show that there holds $J'(v)=0$. In fact, by Lemma \ref{lem.h}, Lemma \ref{lem.embed-3}, Lemma \ref{lem.embed-Lc} and Lebesgue Dominated Convergence Theorem, for any $\psi\in C_0^{\infty}({\mathbb{R}^N})$, we have
\begin{eqnarray*}
&&\langle J'(v_n)-J'(v),\psi\rangle\\
&=&\int_{{\mathbb{R}}^N}(\nabla v_n-\nabla v)\nabla\psi
-\lambda\int_{{\mathbb{R}}^N}c(x)(h(v_n)h'(v_n)-h(v)h'(v))\psi\\
&&-\int_{{\mathbb{R}}^N}(|h(v_n)|^{q-2}h(v_n)h'(v_n)-|h(v)|^{q-2}h(v)h'(v))\psi\\
&&-\int_{{\mathbb{R}}^N}(|h(v_n)|^{2^*-2}h(v_n)h'(v_n)-|h(v)|^{2^*-2}h(v)h'(v))\psi
\to0.
\end{eqnarray*}
Note that $\langle J'(v_n),\psi\rangle\to0$, we get $J'(v)=0$.

In order to prove that $v$ is a weak solution of (\ref{eq.sch5-3}), we must show that $v$ is nontrivial.

\begin{prop}\label{prop.ps-3}
Let $\{v_n\}$ be a Palais-Smale sequence for $J$ at level $c<\frac{1}{N}S^{N/2}$, assume that $v_n\rightharpoonup v$ in $X$, then $v\neq0$.
\end{prop}
{\bf Proof}\quad
We prove the proposition by contradiction. Assume that $v=0$. By Proposition \ref{prop.ps-1}, $\{h(v_n)\}$ is bounded in $X$.

Claim 1: $\{v_n\}$ is also a (PS) sequence for the functional $\tilde{J}: X\to{\mathbb{R}}$ defined by
\begin{eqnarray*}\label{def.J00}
\tilde{J}(v)=\frac{1}{2}\int_{{\mathbb{R}}^N}|\nabla v|^2\mbox{d}x-\frac{1}{q}\int_{{\mathbb{R}}^N}|h(v)|^q\mbox{d}x
-\frac{1}{2^*}\int_{{\mathbb{R}}^N}|h(v)|^{2^*}.
\end{eqnarray*}
Indeed, since the imbedding from $X$ into $L^2({\mathbb{R}^N},c(x))$ is compact, we have
\[
|J(v_n)-\tilde{J}(v_n)|=\frac{\lambda}{2}\int_{{\mathbb{R}}^N}c(x)h(v_n)^2\mbox{d}x\to0,
\]
and for any $\psi\in X$,
\[
|\langle J'(v_n)-\tilde{J}'(v_n),\psi\rangle|=\bigg|\lambda\int_{{\mathbb{R}}^N}c(x)h(v_n)h'(v_n)\psi\bigg|\to0.
\]

Claim 2: For all $R>0$,
\begin{eqnarray}\label{prop.ps-3-2}
\lim\limits_{n\to\infty}\sup\limits_{y\in{\mathbb{R}^N}}\int_{B_R(y)}|h(v_n)|^{q_0}\mbox{d}x=0,
\end{eqnarray}
cannot occur, where $q_0=\max\{2, 2^*(2\alpha)\}$.

Suppose by contradiction that (\ref{prop.ps-3-2}) occurs, that is, $\{v_n\}$ vanished; then by H\"{o}lder's inequality and Sobolev's  inequality,
we have
\begin{eqnarray*}
\|h(v_n)\|_{L^s(B(y,R))}&\leq&\|h(v_n)\|_{L^{q_0}(B(y,R))}^{1-\theta}\|h(v_n)\|_{L^{2^*}(B(y,R))}^{\theta}\\
&\leq& C\|h(v_n)\|_{L^{q_0}(B(y,R))}^{1-\theta}\|\nabla h(v_n)\|_{2}^{\theta},
\end{eqnarray*}
where $\theta=\frac{s-q_0}{2^*-q_0}\frac{2^*}{s}$. Choosing $\theta=2/s$, we obtain
\[
\int_{B(y,R)}|h(v_n)|^s\mbox{d}x\leq C^s\|h(v_n)\|_{L^{q_0}(B(y,R))}^{(1-\theta)s}\int_{B(y,R)}|\nabla h(v_n)|^2\mbox{d}x.
\]
Now covering ${\mathbb{R}}^N$ by balls of radius $R$ in such a way that each point of ${\mathbb{R}}^N$ is contained in at most $N+1$ balls,
we find
\[
\int_{{\mathbb{R}}^N}|h(v_n)|^s\mbox{d}x\leq (N+1)C^s\sup\limits_{y\in{\mathbb{R}}^N}\bigg(\int_{B(y,R)}|h(v_n)|^{q_0}\mbox{d}x\bigg)^{(1-\theta)s/{q_0}}\int_{{\mathbb{R}}^N}|\nabla h(v_n)|^2\mbox{d}x,
\]
which implies that $h(v_n)\to0$ in $L^s({\mathbb{R}}^N)$.
Since $q_0<s<2^*$, by H\"{o}lder's inequality, we get
\begin{eqnarray}\label{prop.ps-3-3}
h(v_n)\to0 ~ in~ L^p({\mathbb{R}}^N),\quad for ~all~ q_0<p<2^*.
\end{eqnarray}
Especially, $h(v_n)\to0$ in $L^q({\mathbb{R}}^N)$.
Now let $\psi=h(v_n)/h'(v_n)$. Since by Lemma \ref{lem.embed-Lc}, $\lambda\int_{\mathbb{R}^N}c(x)h(v_n)^2\mbox{d}x\to0$, we have
\begin{eqnarray*}
o(1)&=&\langle J'(v_n),\psi\rangle\\
&\geq&\int_{{\mathbb{R}}^N}|\nabla h(v_n)|^2\mbox{d}x-\lambda\int_{{\mathbb{R}}^N}c(x)h(v_n)^2\mbox{d}x\\
&&-\int_{{\mathbb{R}}^N}|h(v_n)|^q\mbox{d}x-\int_{{\mathbb{R}}^N}|h(v_n)|^{2^*}\mbox{d}x\\
&\geq&\|h(v_n)\|_X^2-\|h(v_n)\|_{2^*}^{2^*}.
\end{eqnarray*}
By Sobolev's inequality,
\[
o(1)\geq\|h(v_n)\|_{X}^2(1-S^{-2^*/2}\|h(v_n)\|_{X}^{2^*-2}).
\]
If $\|h(v_n)\|_{X}\to0$, then by (5) of Lemma \ref{lem.h}, (\ref{lem.embed-Lc-1}) in Lemma \ref{lem.embed-Lc}, (\ref{prop.ps-3-3}) and Sobolev's inequality,
\begin{eqnarray*}
\int_{{\mathbb{R}}^N}|\nabla v_n|^2\mbox{d}x
&=&\langle J'(v_n),v_n\rangle +\lambda\int_{{\mathbb{R}}^N}c(x)h(v_n)h'(v_n)v_n\mbox{d}x\\
&&+\int_{{\mathbb{R}}^N}|h(v_n)|^{q-2}h(v_n)h'(v_n)v_n\mbox{d}x\\
&&+\int_{{\mathbb{R}}^N}|h(v_n)|^{2^*-2}h(v_n)h'(v_n)v_n\mbox{d}x\\
&\leq&\langle J'(v_n),v_n\rangle +\frac{\lambda}{2\alpha}\int_{{\mathbb{R}}^N}c(x)h(v_n)^2\mbox{d}x\\
&&+\frac{1}{2\alpha}\int_{{\mathbb{R}}^N}|h(v_n)|^{q}\mbox{d}x+\frac{1}{2\alpha}\int_{{\mathbb{R}}^N}|h(v_n)|^{2^*}\mbox{d}x\to0,
\end{eqnarray*}
we contradict $J(v_n)\to c>0$; therefore,
\[
\|h(v_n)\|_{2^*}^{2^*}\geq\|h(v_n)\|_{X}^2+o(1)\geq S^{N/2}+o(1).
\]
By (5) of Lemma \ref{lem.h}, we get
\begin{eqnarray*}
c&=&\lim\limits_{n\to\infty}\bigg\{J(v_n)-\frac{1}{2}\langle J'(v_n),v_n\rangle\bigg\}\\
&=&\lim\limits_{n\to\infty}\bigg\{\frac{\lambda}{2}\int_{{\mathbb{R}}^N}c(x)h(v_n)(h'(v_n)v_n-h(v_n))\mbox{d}x\\
&&+\int_{{\mathbb{R}}^N}|h(v_n)|^{q-2}\Big(\frac{1}{2}h'(v_n)v_n-\frac{1}{q}h(v_n)^2\Big)\mbox{d}x\\
&&+\int_{{\mathbb{R}}^N}|h(v_n)|^{2^*-2}\Big(\frac{1}{2}h'(v_n)v_n-\frac{1}{2^*}h(v_n)^2\Big)\mbox{d}x\bigg\}\\
&\geq&\lim\limits_{n\to\infty}\Big(\frac{1}{2}-\frac{1}{2^*}\Big)\int_{{\mathbb{R}}^N}|h(v_n)|^{2^*}\mbox{d}x\\
&\geq&\frac{1}{N}S^{N/2}
\end{eqnarray*}
which contradicts $c<\frac{1}{N}S^{N/2}$. Thus $\{v_n\}$ does not vanish and there exist $R>0$, $b>0$ and $\{y_n\}\subset{\mathbb{R}^N}$ such that
\[
\lim\limits_{n\to\infty}\int_{B(y,R)}|h(v_n)|^q\mbox{d}x\geq b>0.
\]

Define $\tilde{v}_n(x)=v_n(x+y_n)$. Since $\{v_n\}$ is a (PS) sequence for $\tilde{J}$, $\tilde{v}_n$ is also a (PS) sequence for $\tilde{J}$.
Arguing as in the case of $\{v_n\}$, we get $\tilde{v}_n\rightharpoonup\tilde{v}\in X$ with $\tilde{J}'(\tilde{v})=0$.
Since $\{\tilde{v}_n\}$ does not vanish, we have $\tilde{v}\neq0$. Therefore, by Fatau's lemma, we have
\begin{eqnarray*}
c&\geq&\liminf\limits_{n\to\infty}\bigg\{\tilde{J}(\tilde{v}_n)-\frac{1}{2}\langle\tilde{J}'(\tilde{v}_n),\tilde{v}_n\rangle\bigg\}\\
&=&\liminf\limits_{n\to\infty}\bigg\{\int_{{\mathbb{R}}^N}|h(v_n)|^{q-2}\Big(\frac{1}{2}h'(v_n)v_n-\frac{1}{q}h(v_n)^2\Big)\mbox{d}x\\
&&+\int_{{\mathbb{R}}^N}|h(v_n)|^{2^*-2}\Big(\frac{1}{2}h'(v_n)v_n-\frac{1}{2^*}h(v_n)^2\Big)\mbox{d}x\bigg\}\\
&\geq&\int_{{\mathbb{R}}^N}|h(v_n)|^{q-2}\Big(\frac{1}{2}h'(v)v-\frac{1}{q}h(v)^2\Big)\mbox{d}x\\
&&+\int_{{\mathbb{R}}^N}|h(v_n)|^{2^*-2}\Big(\frac{1}{2}h'(v)v-\frac{1}{2^*}h(v)^2\Big)\mbox{d}x\\
&=&\tilde{J}(\tilde{v})-\frac{1}{2}\langle\tilde{J}'(\tilde{v}),\tilde{v}\rangle.
\end{eqnarray*}
Thus $\tilde{v}\neq0$ is a critical point of $\tilde{J}$ with $\tilde{J}(\tilde{v})\leq c$.

Define
\[
\tilde{c}=\inf_{\gamma\in\tilde{\Gamma}}\sup_{t\in[0,L]}\tilde{J}(\gamma(t))>0,
\]
where $\tilde{\Gamma}=\{\gamma\in C([0,L],X):\gamma(0)=0,\gamma(L)\neq0,\tilde{J}(\gamma(L))<0\}$ for some $L>1$.
Now we construct a path $\gamma(t):[0,L]\to X$ like \cite{JeTa03} such that
\begin{eqnarray*}
\left\{
  \begin{array}{ll}
    \gamma(0)=0, \tilde{J}(\gamma(L))<0, & \hbox{$v\in \gamma([0,L])$;} \\
    \gamma(t)(x)>0, & \hbox{$\forall x\in{\mathbb{R}^N},t\in[0,L]$;} \\
    \max\limits_{t\in[0,L]}\tilde{J}(\gamma(t))=\tilde{J}(\tilde{v}). & \hbox{}
  \end{array}
\right.
\end{eqnarray*}
Setting
\begin{eqnarray}\label{prop.ps-3-4}
\gamma(t)(x)=\left\{
               \begin{array}{ll}
                 \tilde{v}(x/t), & \hbox{$t>0$;} \\
                 0, & \hbox{$t=0$.}
               \end{array}
             \right.
\end{eqnarray}
We see that $\gamma(t)(x)\in\tilde{\Gamma}$ and
\begin{eqnarray*}
&\|\nabla \gamma(t)\|_2^2=t^{N-2}\|\nabla \tilde{v}\|_2^2,\\
&\|h(\gamma(t))\|_q^q=t^{N}\|h(\tilde{v})\|_q^q,\\
&\|h(\gamma(t))\|_{2^*}^{2^*}=t^{N}\|h(\tilde{v})\|_{2^*}^{2^*}.
\end{eqnarray*}
Thus
\[
\tilde{J}(\gamma(t))=\frac{1}{2}t^{N-2}\|\nabla \tilde{v}\|_2^2-t^{N}\Big(\frac{1}{q}\|h(\tilde{v})\|_q^q+\frac{1}{2^*}\|h(\tilde{v})\|_{2^*}^{2^*}\Big).
\]
$\tilde{J}'(\tilde{v})=0$ implies that $\gamma(1)$ is a critical point of $\tilde{J}(\gamma(t))$. Thus $\frac{\mbox{d}}{\mbox{d}t}\Big|_{t=1}\tilde{J}(\gamma(t))=0$. It follows that
\[
\frac{N-2}{2N}\int_{{\mathbb{R}}^N}|\nabla \tilde{v}|^2\mbox{d}x=\frac{1}{q}\int_{{\mathbb{R}}^N}|h(\tilde{v})|^q\mbox{d}x
+\frac{1}{2^*}\int_{{\mathbb{R}}^N}|h(\tilde{v})|^{2^*}\mbox{d}x.
\]
Then
\begin{eqnarray*}
\frac{\mbox{d}}{\mbox{d}t}\tilde{J}(\gamma(t))&=&\frac{N-2}{2}t^{N-3}\int_{{\mathbb{R}}^N}|\nabla \tilde{v}|^2\mbox{d}x\\
&&-Nt^{N-1}\bigg(\frac{1}{q}\int_{{\mathbb{R}}^N}|h(\tilde{v})|^q\mbox{d}x
+\frac{1}{2^*}\int_{{\mathbb{R}}^N}|h(\tilde{v})|^{2^*}\mbox{d}x\bigg)\\
&=&\frac{N-2}{2}t^{N-3}\int_{{\mathbb{R}}^N}|\nabla \tilde{v}|^2\mbox{d}x-\frac{N-2}{2}t^{N-1}\int_{{\mathbb{R}}^N}|\nabla \tilde{v}|^2\mbox{d}x\\
&=&\frac{N-2}{2}t^{N-3}(1-t^2)\int_{{\mathbb{R}}^N}|\nabla \tilde{v}|^2\mbox{d}x.
\end{eqnarray*}
We conclude that $\frac{\mbox{d}}{\mbox{d}t}\tilde{J}(\gamma(t))>0$ for $t\in(0,1)$ and $\frac{\mbox{d}}{\mbox{d}t}\tilde{J}(\gamma(t))<0$ for $t\in(1,L)$. Thus we get the desired path.

If $\lambda=0$, we have proved the proposition. For $\lambda>0$, since the path $\gamma$ given by (\ref{prop.ps-3-4}) belongs to $\tilde{\Gamma}\subset\Gamma$ after scaling, we obtain
\[
c\leq\max\limits_{t\in[0,L]}J(\gamma(t))=J(\gamma(\overline{t}))<\tilde{J}(\gamma(\overline{t}))
\leq\max\limits_{t\in[0,L]}\tilde{J}(\gamma(t))=\tilde{J}(\tilde{v})\leq c,
\]
which is a contradiction. Therefore, $v$ is nontrivial.
$\quad\Box$

\section{Proof of main theorems}
\setcounter{equation}{0}
\label{sec:5-4}

In this section, we will study the properties of the functional $J$ and prove the main theorem, this include the construction of a path that has a maximum level $c<\frac{1}{N}S^{N/2}$.

\begin{lem}\label{lem.h-2*-1}
There exists $d_0>0$ such that
\[
\lim\limits_{t\to+\infty}(t-h(t))\geq d_0.
\]
\end{lem}
{\bf Proof}\quad
Assume that $t>0$. By Lemma \ref{lem.h} we have $h(t)\leq t$ and $h(t)\leq h'(t)t$. Thus we have
\begin{eqnarray*}
t-h(t)&\geq&t(1-h'(t))\\
&=&t\frac{(2\kappa\alpha^2+h(t)^{2(1-2\alpha)})^{1/2}-h(t)^{1-2\alpha}}
{(2\kappa\alpha^2+h(t)^{2(1-2\alpha)})^{1/2}}\\
&\geq&\frac{\kappa\alpha^2t}{2\kappa\alpha^2+h(t)^{2(1-2\alpha)}}\\
&\geq&\frac{\kappa\alpha^2t}{2h(t)^{2(1-2\alpha)}}\quad\mbox{for $t$ large}\\
&:=&d(\alpha,t).
\end{eqnarray*}

Case 1. If $\frac{1}{4}<\alpha<\frac{1}{2}$, then $0<1-2\alpha<\frac{1}{2}$ and thus $d(\alpha,t)\to+\infty$ as $t\to+\infty$.

Case 2. If $\alpha=\frac{1}{4}$, then $1-2\alpha=1$ and thus $d(\alpha,t)\to\frac{\kappa\alpha^2}{2}$ as $t\to+\infty$.

Case 3. If $0<\alpha<\frac{1}{4}$, we claim that $t-h(t)\to0$ is impossible.
Assume on the contrary. Note that $4\alpha<1$ and $h(t)^{4\alpha-1}\to0$ as $t\to+\infty$, by L'Hospital's Principle, we have
\begin{eqnarray*}
0&\leq&\lim\limits_{t\to+\infty}\frac{t-h(t)}{h(t)^{4\alpha-1}}\\
&=&\lim\limits_{t\to+\infty}\frac{1-h'(t)}{(4\alpha-1)h(t)^{4\alpha-2}h'(t)}\\
&=&\lim\limits_{t\to+\infty}\frac{h(t)^{1-2\alpha}}{4\alpha-1}
[(2\kappa\alpha^2+h(t)^{2(1-2\alpha)})^{1/2}-h(t)^{1-2\alpha}]\\
&=&\frac{\kappa\alpha^2}{4\alpha-1}<0,
\end{eqnarray*}
a contradiction.

In summation, for all $0<\alpha<1/2$, there exists $d_0>0$ such that the conclusion of the lemma holds.
$\quad\Box$

\begin{lem}\label{lem.h-2*-2}
For $h(t)$ defined in (\ref{def.h}), we have

(i) If $\frac{1}{4}<\alpha<\frac{1}{2}$, then
\[
\lim\limits_{t\to+\infty}\frac{t-h(t)}{t^{4\alpha-1}}=\frac{\kappa\alpha^2}{4\alpha-1};
\]

(ii) If $0<\alpha\leq\frac{1}{4}$, then
\[
\lim\limits_{t\to+\infty}\frac{t-h(t)}{\log h(t)}
\leq\left\{
   \begin{array}{ll}
     \frac{\kappa}{16}, & \hbox{$\alpha=\frac{1}{4}$;} \\
     0, & \hbox{$0<\alpha<\frac{1}{4}$.}
   \end{array}
 \right.
\]
\end{lem}
{\bf Proof}\quad
(i) Assume that $t>0$. By the proof of Lemma \ref{lem.h-2*-1}, when $\frac{1}{4}<\alpha<\frac{1}{2}$, we have $t-h(t)\to+\infty$ as $t\to+\infty$, then we can use L'Hospital's Principle to compute that
\begin{eqnarray*}
\lim\limits_{t\to+\infty}\frac{t-h(t)}{t^{4\alpha-1}}
=\lim\limits_{t\to+\infty}\frac{1-h'(t)}{(4\alpha-1)t^{4\alpha-2}}
=\frac{\kappa\alpha^2}{4\alpha-1}
\end{eqnarray*}
(ii) When $0<\alpha\leq\frac{1}{4}$ and if there exists a constant $C>0$ such that $t-h(t)\leq C$, then the conclusion holds. Otherwise, we may assume that $t-h(t)\to+\infty$ as $t\to+\infty$. Then again by L'Hospital's Principle, we have
\begin{eqnarray*}
A&:=&\lim\limits_{t\to+\infty}\frac{t-h(t)}{\log h(t)}\\
&=&\lim\limits_{t\to+\infty}h(t)\bigg(\frac{1}{h'(t)}-1\bigg)\\
&=&\lim\limits_{t\to+\infty}\frac{2\kappa\alpha^2h(t)^{2\alpha}}
{(2\kappa\alpha^2+h(t)^{2(1-2\alpha)})^{1/2}+h(t)^{1-2\alpha}}.
\end{eqnarray*}
Thus $A=\frac{\kappa}{16}$ when $\alpha=\frac{1}{4}$ and $A=0$ when $0<\alpha<\frac{1}{4}$.
This completes the proof.
$\quad\Box$

To finish the proof of Theorem \ref{thm.m1}, we construct a path which minimax level less than $\frac{1}{N}S^{N/2}$.

\begin{prop}\label{prop.mp-1}
The minimax level $c$ defined in (\ref{def.c}) satisfies
\[
c<\frac{1}{N}S^{N/2}.
\]
\end{prop}
{\bf Proof}\quad
We follow the strategy used in \cite{BrNi93}. Let
\[
v^*=\frac{[N(N-2)\varepsilon^2]^{(N-2)/4}}{(\varepsilon^2+|x|^2)^{(N-2)/2}}
\]
be the solution of $-\Delta u=u^{2^*-1}$ in $\mathbb{R}^N$. Then
\begin{eqnarray*}
\int_{{\mathbb{R}}^N}|\nabla v^*|^2=\int_{{\mathbb{R}}^N}|v^*|^{2^*}=S^{N/2},
\end{eqnarray*}
Let $\eta_\varepsilon(x)\in C_0^\infty({\mathbb{R}^N},[0,1])$ be a cut-off function with
$\eta_\varepsilon(x)=1$ in $B_\varepsilon=\{x\in{\mathbb{R}^N}:|x|\leq\varepsilon\}$ and
$\eta_\varepsilon(x)=0$ in $B^c_{2\varepsilon}={\mathbb{R}^N}\setminus B_{2\varepsilon}$.
Let $v_\varepsilon=\eta_\varepsilon v^*$.
For all $\varepsilon>0$, there exists $t^\varepsilon>0$ such that $J(t^\varepsilon v_\varepsilon)<0$ for all $t>t^\varepsilon$. Define the class of paths
\begin{eqnarray*}
\Gamma_\varepsilon=\{\gamma\in C([0,1],X):\gamma(0)=0,\gamma(1)=t^\varepsilon v_\varepsilon\}
\end{eqnarray*}
and the minimax level
\begin{eqnarray*}
c_\varepsilon=\inf_{\gamma\in\Gamma_\varepsilon}\max_{t\in[0,1]}J(\gamma(t))
\end{eqnarray*}
Let $t_\varepsilon$ be such that
\[
J(t_\varepsilon v_\varepsilon)=\max_{t\geq0}J(tv_\varepsilon)
\]
Note that the sequence $\{v_\varepsilon\}$ is uniformly bounded in $X$, we conclude that $\{t_\varepsilon\}$ is upper and lower bounded by two positive constants. In fact, if $t_\varepsilon\to0$, we have $J(t_\varepsilon v_\varepsilon)\to0$; otherwise, if $t_\varepsilon\to+\infty$, we have $J(t_\varepsilon v_\varepsilon)\to-\infty$. In both cases we get contradictions according to Lemma \ref{lem.mp-1}. This proved the conclusion.

According to \cite{BrNi93}, we have, as $\varepsilon\to0$,
\begin{eqnarray}\label{prop.v-1}
\|\nabla v_\varepsilon\|^2_{2}=S^{N/2}+O(\varepsilon^{N-2}),\quad
\|v_\varepsilon\|^{2^*}_{{2^*}}=S^{N/2}+O(\varepsilon^N).
\end{eqnarray}

Define
\begin{eqnarray*}
H(t_\varepsilon v_\varepsilon)
=-\frac{\lambda}{2}\int_{\mathbb{R}^N}c(x)h(t_\varepsilon v_\varepsilon)
-\frac{1}{q}\int_{\mathbb{R}^N}h(t_\varepsilon v_\varepsilon)^q
+\frac{1}{2^*}\int_{\mathbb{R}^N}[(t_\varepsilon v_\varepsilon)^{2^*}
-h(t_\varepsilon v_\varepsilon)^{2^*}]
\end{eqnarray*}
By the definition of $v_\varepsilon$, for $x\in B_\varepsilon$, there exist two constants $c_2\geq c_1>0$ such that for $\varepsilon$ small enough, we have
\[
c_1\varepsilon^{-(N-2)/2}\leq v_\varepsilon(x)\leq c_2\varepsilon^{-(N-2)/2}
\]
and
\[
c_1\varepsilon^{-(N-2)/2}\leq h(v_\varepsilon(x))\leq c_2\varepsilon^{-(N-2)/2}.
\]
Note that $t_\varepsilon$ is upper and lower bounded, $c(x)$ is continuous in $\overline{B}_\varepsilon$, there exist constants $C_1>0, C_2>0$ such that
\begin{eqnarray}\label{prop.mp-1-2}
\int_{B_\varepsilon}c(x)h^2(t_\varepsilon v_\varepsilon)\geq C_1\varepsilon^{2} =C_1\varepsilon^{(\frac{2^*}{2}-1)(N-2)}
\end{eqnarray}
and
\begin{eqnarray}\label{prop.mp-1-3}
\int_{B_\varepsilon}h^q(t_\varepsilon v_\varepsilon)\geq C_2\varepsilon^{N-q\frac{N-2}{2}} =C_2\varepsilon^{(\frac{2^*}{2}-\frac{q}{2})(N-2)}.
\end{eqnarray}
Moreover, note that $h(t_\varepsilon v_\varepsilon)\leq t_\varepsilon v_\varepsilon$ and $2^*>2$, by H\"{o}lder's inequality, we have
\begin{eqnarray*}
R_\varepsilon&:=&
\frac{1}{2^*}\int_{B_\varepsilon}[(t_\varepsilon v_\varepsilon)^{2^*}-h^{2^*}(t_\varepsilon v_\varepsilon)]\\
&\leq& \int_{B_\varepsilon}(t_\varepsilon v_\varepsilon)^{2^*-1}(t_\varepsilon v_\varepsilon-h(t_\varepsilon v_\varepsilon))\\
&\leq&\bigg(\int_{B_\varepsilon}(t_\varepsilon v_\varepsilon)^{2^*}\bigg)^{\frac{2^*-1}{2^*}} \bigg(\int_{B_\varepsilon}(t_\varepsilon v_\varepsilon-h(t_\varepsilon v_\varepsilon))^{2^*}\bigg)^{\frac{1}{2^*}}.
\end{eqnarray*}
According to Lemma \ref{lem.h-2*-2}, there exist $C_3>0$ such that for $\frac{1}{4}<\alpha<\frac{1}{2}$,
\begin{eqnarray}\label{prop.mp-1-4}
R_\varepsilon\leq C_3\bigg(\int_{B_\varepsilon}(t_\varepsilon v_\varepsilon)^{2^*(4\alpha-1)}\bigg)^{\frac{1}{2^*}}
\leq C_3\varepsilon^{(1-2\alpha)(N-2)},
\end{eqnarray}
while for $0<\alpha<\frac{1}{4}$, there exists a constant $\delta\in(0,1)$ such that
\begin{eqnarray}\label{prop.mp-1-5}
R_\varepsilon\leq C_3\bigg(\int_{B_\varepsilon}(t_\varepsilon v_\varepsilon)^{2^*\delta}\bigg)^{\frac{1}{2^*}}
\leq C_3\varepsilon^{\frac{1}{2}(1-\delta)(N-2)}.
\end{eqnarray}
From the above estimations (\ref{prop.mp-1-2})-(\ref{prop.mp-1-5}), we get
\begin{eqnarray}\label{prop.mp-1-6}
H(t_\varepsilon v_\varepsilon)\leq -C_1\varepsilon^{(\frac{2^*}{2}-1)(N-2)} -C_2\varepsilon^{(\frac{2^*}{2}-\frac{q}{2})(N-2)} +C_3\varepsilon^{(1-2\alpha)(N-2)}
\end{eqnarray}
when $\frac{1}{4}<\alpha<\frac{1}{2}$ and
\begin{eqnarray}\label{prop.mp-1-7}
H(t_\varepsilon v_\varepsilon)\leq -C_1\varepsilon^{(\frac{2^*}{2}-1)(N-2)} -C_2\varepsilon^{(\frac{2^*}{2}-\frac{q}{2})(N-2)}
+C_3\varepsilon^{\frac{1}{2}(1-\delta)(N-2)}
\end{eqnarray}
when $0<\alpha<\frac{1}{4}$.

Now we have
\begin{eqnarray}\label{prop.mp-1-1}
J(t_\varepsilon v_\varepsilon)=\frac{t_\varepsilon^2}{2}\int_{\mathbb{R}^N}|\nabla v_\varepsilon|^2-\frac{t_\varepsilon^{2^*}}{2^*}\int_{\mathbb{R}^N}|v_\varepsilon|^{2^*} +H(t_\varepsilon v_\varepsilon).
\end{eqnarray}
Since the function $\xi(t)=\frac{1}{2}t^2-\frac{1}{2^*}t^{2^*}$ achieve its maximum $\frac{1}{N}$ at point $t_0=1$, by using (\ref{prop.v-1}), we derive from (\ref{prop.mp-1-1}) that
\begin{eqnarray}\label{prop.mp-1-8}
J(t_\varepsilon v_\varepsilon)\leq \frac{1}{N}S^{N/2}+H(t_\varepsilon v_\varepsilon)+O(\varepsilon^{N-2}).
\end{eqnarray}
Combining (\ref{prop.mp-1-6}), (\ref{prop.mp-1-7}) and (\ref{prop.mp-1-8}), we conclude that

(i) for $\frac{1}{4}<\alpha<\frac{1}{2}$ and $q>\frac{4}{N-2}+4\alpha$, we have $(\frac{2^*}{2}-\frac{q}{2})(N-2)<(1-2\alpha)(N-2)$;

(ii) for $0<\alpha<\frac{1}{4}$ and $q>\frac{N+2}{N-2}$, we have $(\frac{2^*}{2}-\frac{q}{2})(N-2)<\frac{1}{2}(1-\delta)(N-2)$ for $\delta>0$ small enough.

Therefore, conclusions (i)-(ii) give that
\begin{eqnarray}\label{prop.mp-1-9}
c_\varepsilon=J(t_\varepsilon v_\varepsilon)< \frac{1}{N}S^{N/2}.
\end{eqnarray}
Finally, since $\Gamma_\varepsilon\subset\Gamma$, we have
\[
c\leq c_\varepsilon<\frac{1}{N}S^{N/2}.
\]
This proved the proposition.
$\quad\Box$

{\bf Proof of Theorem \ref{thm.m1}}\quad
Firstly, by Lemma \ref{lem.mp-1}-\ref{lem.mp-2}, the functional $J$ has the Mountain Pass Geometry. Then there exists a Palais-Smale sequence $\{v_n\}$ at level $c$ given in (\ref{def.c}). Secondly, by Proposition \ref{prop.ps-1}, the Palais-Smale sequence $\{v_n\}$ is bounded in $X$. By Proposition \ref{prop.ps-3}, if $c<\frac{1}{N}S^{N/2}$, then the weak limit $v$ of $\{v_n\}$ in $X$ is nonzero and it is a critical point of $J$. Finally, by Proposition \ref{prop.mp-1}, there indeed exists a mountain pass which maximum level $c_\varepsilon$ is strictly less than $\frac{1}{N}S^{N/2}$. This implies that the level $c<\frac{1}{N}S^{N/2}$ and $v$ is a nontrivial weak solution of Eq.(\ref{eq.sch5-h}). Then $u=h(v)$ is a weak solution of Eq.(\ref{eq.sch5-3}).
$\quad\Box$

{\bf Proof of Corollary \ref{rem.m1-3}}\quad
According to (\ref{lem.embed-2-0}), we have
\[
|\nabla (|h(v)|^a)|^2\leq\frac{a^2|h(v)|^{2(a-2\alpha)}}{2\alpha+|h(v)|^{2(1-2\alpha)}}|\nabla v|^2.
\]
Let
\[
\eta(s)=\frac{s^{2(a-2\alpha)}}{2\alpha+s^{2(1-2\alpha)}},\quad s\in(0,+\infty),
\]
then $\eta(s)$ has a unique maximum point $s_0$ satisfies $s_0^{2(1-2\alpha)}=\frac{2\alpha(a-2\alpha)}{(1-a)}$.
Assume that $s_0\geq1$, that is, $a_0=\frac{1+(2\alpha)^2}{1+2\alpha}\leq a\leq1$, then we have
\[
\eta(s_0)=\frac{s_0^{2(a-2\alpha)}}{2\alpha+s_0^{2(1-2\alpha)}}\leq\frac{s_0^{2(1-2\alpha)}}{2\alpha+s_0^{2(1-2\alpha)}}\leq1.
\]
This implies $|\nabla (|h(v)|^a)|^2\leq|\nabla v|^2$, thus the constant $C_1=1$ in (\ref{lem.mp-1-1}). Then conclusion follows from Theorem \ref{thm.m1}.
$\quad\Box$






\begin{thebibliography}{00}

\bibitem{LiWa03} J. Liu, Z.-Q. Wang, Soliton solutions for quasilinear Schr\"{o}dinger equations, I, {\it Proc. Amer. Math. Soc.},
{\bf 131} (2003) 441-448.

\bibitem{LiWW03} J. Liu, Y. Wang, Z.-Q. Wang, Soliton solutions for quasilinear Schr\"{o}dinger equations, II, {\it J. Differential Equations}, {\bf 187} (2003) 473-493.


\bibitem{Kuri81} S. Kurihura, Large-amplitude quasi-solitons in superfluid films, {\it J. Phys. Soc. Japan}, {\bf 50} (1981) 3262-3267.

\bibitem{LiSe78} A. G. Litvak, A. M. Sergeev, One dimensional collapse of plasma waves, {\it JETP Lett.}, {\bf 27} (1978) 517-520.

\bibitem{KoIK90} A. M. Kosevich, B. A. Ivanov, A.S. Kovalev, Magnetic solitons, {\it Phys. Rep.}, {\bf 194} (1990) 117-238.

\bibitem{QuCa82} G. R. W. Quispel, H. W. Capel, Equation of motion for the Heisenberg spin chain, {\it Physica A}, {\bf 110} (1982) 41-80.

\bibitem{Hass80} R. W. Hasse, A general method for the solution of nonlinear soliton andkink Schr\"{o}dinger equations, {\it Z. Phys.}, {\bf 37} (1980) 83-87.

\bibitem{MaFe84} V. G. Makhankov, V. K. Fedyanin, Non-linear effects in quasi-one-dimensional models of condensed matter theory, {\it Phys. Rep.}, {\bf104} (1984) 1-86.


\bibitem{BoGa93} A. V. Borovskii, A. L. Galkin, Dynamical modulation of an ultrashort high-intensity laser pulse in matter, {\it JETP}, {\bf 77} (1993) 562-573.

\bibitem{CoJe04} M. Colin, L. Jeanjean, Solutions for a quasilinear Schr?dinger equations: a dual approach, {\it Nonlinear Anal. TMA}, {\bf 56} (2004) 213-226.

\bibitem{OMS07} J. M. B. do \'{O}, O H. Miyagaki, S. H. M. Soares, Soliton solutions for quasilinear Schr\"{o}dinger equations: the critical exponential case, {\it Nonlinear Anal. TMA}, {\bf 67} (2007) 3357-3372.

\bibitem{OMS10} J. M. B. do \'{O}, O H. Miyagaki, S. H. M. Soares, Soliton solutions for quasilinear Schr\"{o}dinger equations with critical growth, {\it J. Differential Equations, } {\bf 1248} (2010)   722--744.


\bibitem{SiVi10}E. A. B. Silva, G. F. Vieira, Quasilinear asymptotically periodic Schr\"{o}dinger equations with critical growth, {\it Calc. Var. PDEs,} {\bf 39} (2010) 1-33.

\bibitem{Moam06} A. Moameni, Existence of soliton solutions for a quasilinear Schr\"{o}dinger equation involving critical exponent in $\mathbb{R}^N$, {\it J. Differential Equations,} {\bf 229} (2006) 570-587.


\bibitem{LiZh13} Z. Li, Y. Zhang, Solutions for a class of quasilinear Schrodinger equations with critical exponents term, {\it J. Math. Phys.} {\bf 58} 021501 (2017).
    

\bibitem{BoGa93} A. V. Borovskii, A. L. Galkin, Dynamical modulation of an ultrashort high-intensity laser pulse in matter, {\it JETP}, {\bf 77} (1993) 562-573.

\bibitem{JeTa03} L. Jeanjean, K. Tanaka, A remark on least energy solutions in ${\mathbb{R}^N}$, {\it Proc. Amer. Math. Sco.}, {\bf 131} (2003) 2399-2408.


\bibitem{BrNi93} H. Brezis, L. Nirenberg, Positive solutions of nonlinear elliptic equations involving critical Sobolev exponent, {\it Comm. Pure Appl. Math.}, {\bf 36} (1993) 437-477.





\end{thebibliography}



\end{document}